\documentclass[11pt, reqno]{amsart}
\usepackage{hyperref}

\usepackage[sort, numbers]{natbib}  

 \usepackage{amsthm} 
 \usepackage{amssymb}
 \usepackage{amsmath}
 \usepackage{amsfonts}
 \usepackage{stmaryrd}
 \usepackage{fancyhdr}
 \usepackage{graphicx}
 \usepackage{adjustbox}

 \usepackage{graphicx,color}
 \usepackage{epstopdf}
\usepackage{subfigure}

\newcommand{\N}{\mathbb{N}}
\newcommand{\E}{\mathbb{E}}
\newcommand{\cL}{\mathcal{L}}

\newcommand{\pa}{\partial}

\newcommand{\be}{\begin{equation}}
	\newcommand{\ee}{\end{equation}}

\newcommand{\R}{\mathbb{R}}
\newcommand{\T}{\mathbb T}
\newcommand{\iin}{^{i,N}}
\newcommand{\jn}{^{j,N}}

\newcommand{\Dt}{\Delta t}
\DeclareMathOperator{\diag}{diag}

\newtheorem{theorem}{Theorem}[section]
\newtheorem{note}[theorem]{Note}

\begin{document}

\title[Counting number of solutions of PDEs via infinite dimensional sampling]{Counting the number of stationary solutions of Partial Differential Equations via infinite dimensional sampling}
    \author{M. Kolodziejczyk$^{1}$}
    \address{$^{1}$Mathematics Department, Politecnico di Milano, Milan, Italy, martin.kolodziejczyk@polimi.it}

    \author{M. Ottobre$^{2}$}
    \address{$^{2}$ Corresponding Author. Mathematics Department, Heriot-Watt University and Maxwell Institute, Edinburgh, Scotland, m.ottobre@hw.ac.uk}
\author{G. Simpson$^{3}$}  
\address{$^{3}$Drexel University, Philadelphia, United States, grs53@drexel.edu}


\maketitle

\begin{abstract}
This paper is concerned with the problem of  counting solutions of stationary nonlinear Partial Differential Equations (PDEs) when the PDE is known to admit more than one solution. We suggest tackling the problem via a sampling-based approach. We test our proposed methodology on the McKean-Vlasov PDE, more precisely on the problem of determining the number of stationary solutions of the McKean-Vlasov (or porous medium) equation.


\end{abstract}

\section{Introduction}
 Within the field of Partial Differential Equations (PDEs), a large body of literature is concerned with investigating well-posedness of  PDEs. That is, with studying existence and {\em uniqueness} of solutions, in appropriate spaces.  It is however the case that  many nonlinear PDEs of pivotal importance in applications admit more than one solution. This paper is devoted to the problem of determining the number of solutions  of stationary nonlinear PDEs (and then possibly finding properties of such solutions), i.e. solving problems of the form 
\begin{equation}\label{eqn: stationary PDE}
    \cL u =0 \,, 
\end{equation}
where $\cL$ is a suitable (nonlinear) operator acting on a given function space. Studying the set of solutions of \eqref{eqn: stationary PDE} is central in at least three broad contexts. 
First, solutions of \eqref{eqn: stationary PDE} are candidate equilibria of the associated time-dependent PDE,
\begin{equation}\label{eqn: time-dependent PDE}
    \pa_t u = \cL u \,.
\end{equation}
This is relevant  for both  deterministic and  stochastic problems. Indeed, when $\cL$ is the Fokker-Planck operator associated with a given Stochastic (ordinary) Differential Equation (SDE), 
the solutions of \eqref{eqn: stationary PDE} are the invariant 

\noindent
measures of the SDE. Second, at least when the solutions of \eqref{eqn: stationary PDE} are stable -- in the sense that they are stable equilibria of \eqref{eqn: time-dependent PDE} -- one expects them to correspond to  metastable states of  the associated SPDE, i.e. of the evolution 
\begin{equation}
\label{e:noise1}
    \pa_t u = \cL u + \epsilon \,  \eta(t,x) \, ,
\end{equation}
where $\epsilon>0$ is a small parameter and $\eta$ is spatiotemporal noise.   Finally, if the problem  has a  variational structure, solutions of \eqref{eqn: stationary PDE} are related  to extrema of associated functionals, often referred to as energy functionals or potentials.  Hence the  connection with the many applications of interest in calculus of variations and with energy landscape exploration problems. 

Many problems in both natural and applied sciences can be recast as one of \eqref{eqn: stationary PDE}, \eqref{eqn: time-dependent PDE}, and \eqref{e:noise1}. 
 Let us give some examples.  Stochastic processes with multiple invariant measures naturally arise  in the study of interacting particle systems and associated mean field limits. Such systems are paradigmatic models in statistical mechanics and kinetic theory, in connection with the study of phase transitions, as well as  in  the study of  collective navigation and  consensus formation \cite{talay1996probabilistic, dawson1983critical}. The importance of understanding processes with multiple invariant measures cannot be understated. Indeed,  while statistical sampling has been one of the main  motivations for  the study of processes with one invariant measure -  {\em ergodic processes} - many processes in natural and engineered systems are likely to exhibit multiple invariant measures.  This includes the tendency of fireflies to synchronise their flashing or not, the microstructure of nematic crystals aligning to several distinct equilibrium configurations (giving rise to different material properties), the opinion formation process in social media, which  can converge towards various possible outcomes, etc \cite{goddard2023study, freidlin2021averaging,pavliotis, pareschi2013interacting, gorbonos2024geometrical, yin2001limit}. 

Related to problems of the form \eqref{e:noise1}, metastability theory and landscape exploration are central to molecular dynamics, computational chemistry and condensed matter physics, particularly for the calculation of reaction rates. Indeed, a large number of methods have been proposed to understand transitions between metastable states, including the string method in finite \cite{weinan2002string} and infinite dimensions \cite{weinan2004minimum}, the nudged elastic band method \cite{henkelman2000climbing}, surface walking approaches such as the dimer method, the gentlest ascent dynamics \cite{weinan2011gentlest}, the activation-relaxation technique, and more. All such methods require some prior knowledge of the metastable states (so that one can e.g. `initialise the string'),  which are a subset of the solutions of stationary problems like \eqref{eqn: stationary PDE}.

Despite the pervasiveness of PDEs with multiple solutions, the methods available to find all such solutions is very limited. To the best of our knowledge, aside from the  brute force approach of  starting the Newton algorithm from a large number of  different initial points, there are few systematic strategies that work in full generality. The available tools include deflation \cite{farrell2015deflation} (which only requires pre-knowledge of one solution {and can find any stationary solution, irrespective of its stability}), continuation and homotopy-continuation \cite{mehta2011finding} (which have similar requirements).  These approaches are all numerical or computational as, except for rare cases, one does not expect to be able to find solutions analytically. \footnote{We note in passing that in differential geometry the Atiyah-Singer index theorem \cite{atiyah1970global} serves precisely the purpose of counting solutions of differential equations associated to linear elliptic operators. It might be worth exporting this theorem to computational stochastic analysis and reconciling that framework with what is known in more applied fields. However, such a theorem does not apply to considered in this paper.}  We mention also so-called  `landscape exploration methods' such as the  eigenvector-following method \cite{doye1997surveying} and minimax approaches \cite{li2001minimax}  to find saddle points, which are also relevant to the study of reaction rates. However, these approaches mostly apply to finite dimensional landscapes.

We emphasize that the complexity of the problems that arise in this context is  high; examples are illustrated in the works \cite{wang2019order, han2021solution} concerning the study of nematic crystals (different solutions of the associated stationary PDE problem here correspond to  distinct equilibrium configurations of the crystal),  in the setting of de Gennes theory. In \cite{wang2019order},  twenty eight solutions were found for the problem considered there, with no guarantee that all had been found.   For semilinear Schr\"odinger equations, there are (at least) countably many soliton type solutions, even after accounting for symmetries, \cite{sulem_nonlinear_1999,fibich_nonlinear_2014}.

There are two challenges to determining the set of solutions to \eqref{eqn: stationary PDE}: determining the number of solutions and then finding the solutions themselves (or at least  some properties of such solutions). These challenges should be regarded as distinct.  Deflation is a very powerful method,  proven useful in many settings. But even by deflation one cannot guarantee all solutions have been found. This limitation is true of all the other methods we have mentioned, and it will certainly be true, also, of the method we propose in this paper. Thus, it is pragmatic to have at one's disposal an array of methods and to exploit synergies between them,  to make progress on the specific problem at hand.  For example,  if by deflation one has found $n$ solutions but some `solution counting method' reveals that we should expect at least say $n+2$ then one can redouble the search. {Knowledge of such solutions can, in turn, provide the basis to use paradigms such as e.g. the string method, which requires pre-knowledge of the metastable states (to initialise the string).}  

This paper attempts to put forward an approach to address the problem of counting the number of (stable) solutions of a given stationary PDE of the form \eqref{eqn: stationary PDE}.   The purpose of this short paper is to present the main  idea of this approach  though it is undoubtedly the case that further investigation and (non-trivial) analysis will be needed in order to understand its strengths and limitations. We comment on this in Section \ref{subsec: approach intro} below, Section \ref{sec: background} and at various points throughout.   

The paper is organised as follows: in Subsection \ref{subsec: approach intro} we present our method and then, in subsequent sections, we apply this method to find the number of stable stationary solutions of the McKean-Vlasov (or porous media) equation.  Such stationary solutions correspond to invariant measures of stochastic processes which are non-linear in the sense of McKean. We do this both for parameter regimes where the number of solutions is, a priori, known and also in settings where the issue of finding the  number of solutions is yet to be settled.  In Section \ref{sec: background}  we provide background material on the McKean-Vlasov PDE and  on the associated  McKean-Vlasov SPDE;  the latter will be the central tool we employ  to find the number of stable solutions of the stationary McKean-Vlasov PDE. {Section \ref{sec: background} also serves the purpose of illustrating how determining the number of stationary solutions of a given PDE can be rather non-trivial, even in simple situations}. Section \ref{sec:numerical discretization} and Section \ref{sec:simulations} are devoted to numerical and simulation aspects, respectively (more detail below).

\subsection{Our Approach.} \label{subsec: approach intro} We propose that in order to count the number of stable solutions of the stationary problem \eqref{eqn: stationary PDE} one can sample 
from the  invariant measure (assuming it exists and is unique) of  the Stochastic Partial Differential Equation (SPDE) obtained by  perturbing the time-dependent problem \eqref{eqn: time-dependent PDE} by additive noise; that is, one can sample from the invariant measure of  the SPDE
\begin{equation}\label{eqn: additive SPDE intro}
    \pa_t u = \cL u + \eta(t,x) \, ,
\end{equation}
where $\eta$ is a suitable space-time noise. At least in the case in which the solutions of \eqref{eqn: stationary PDE} are isolated, we expect the number of modes of the invariant measure of \eqref{eqn: additive SPDE intro} to coincide with the number of stable solutions of \eqref{eqn: stationary PDE} and the number of stable equilibria of \eqref{eqn: time-dependent PDE}. Further, we expect a correspondence between the metastable states of \eqref{eqn: additive SPDE intro}, the modes,  and the stable stationary solutions of \eqref{eqn: stationary PDE} and \eqref{eqn: time-dependent PDE}. The noise $\eta$ must be appropriately chosen not only to make sure that \eqref{eqn: additive SPDE intro} is well posed and has a unique invariant measure, but also to, potentially, satisfy other constraints specific to the problem under consideration, see Section \ref{sec: background} for more detail on this point.\footnote{We consider additive noise because it is straightforward from both an analytical and computational perspective, but the use of other types of noise merits exploration.}  Intuitively,  if the noise $\eta$ is strong enough we would expect \eqref{eqn: additive SPDE intro} to have at most one invariant measure. 

\begin{figure}
    \centering
    \subfigure[]{\includegraphics[width=6.5cm]{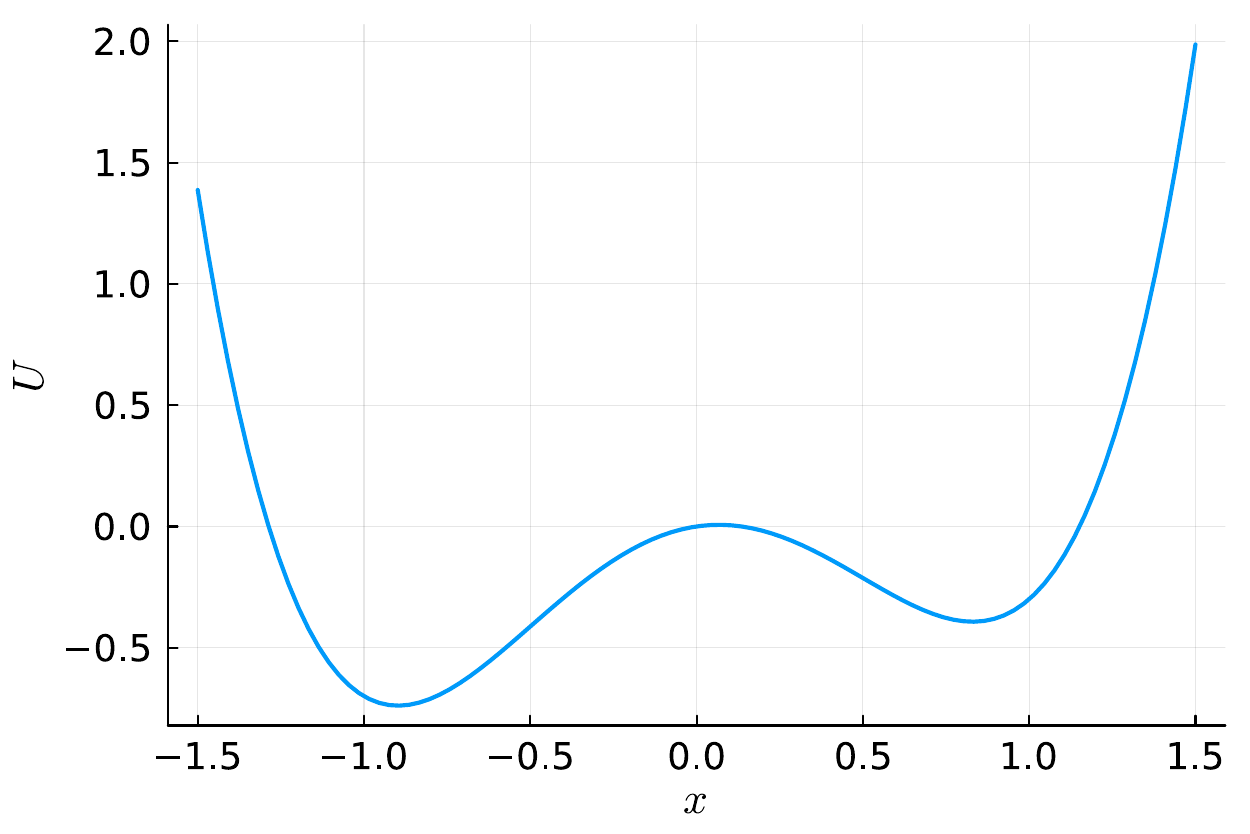}}
    \subfigure[]{\includegraphics[width=6.5cm]{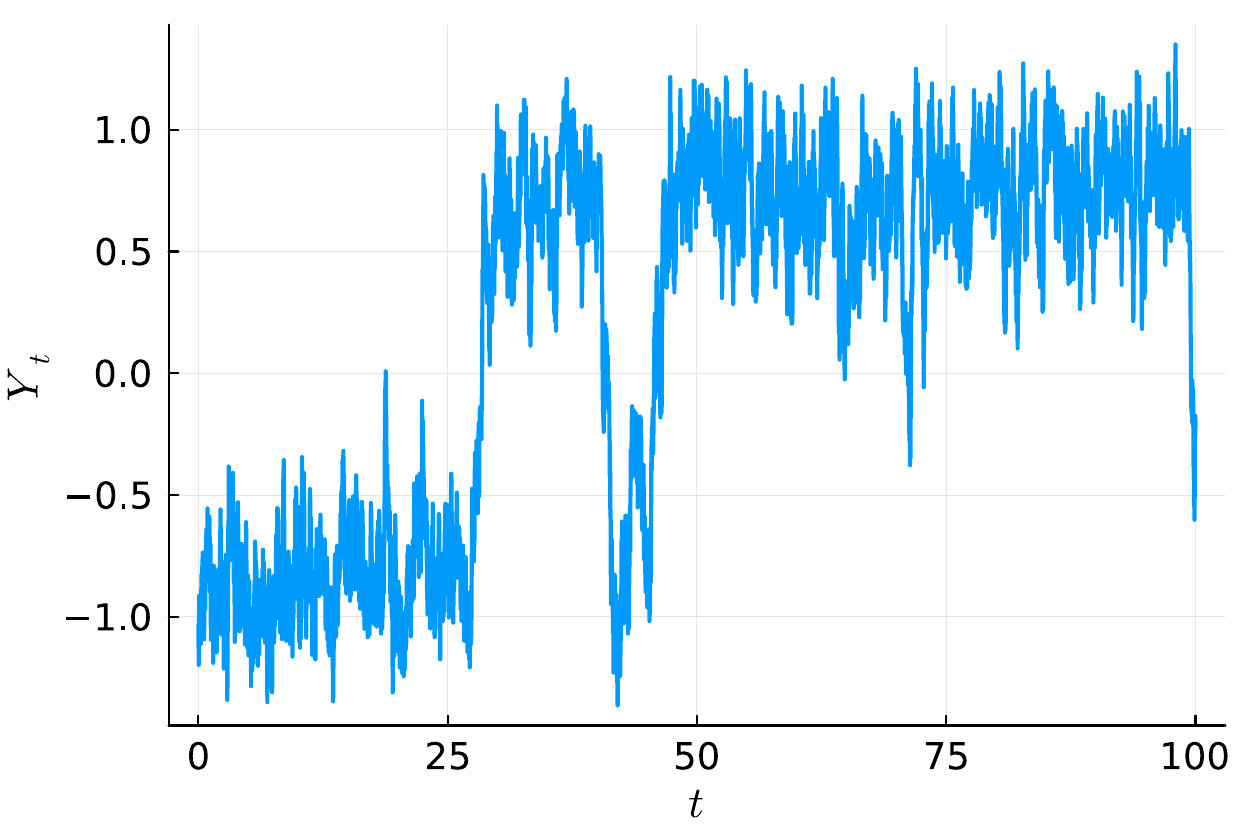}}
    \caption{A simple scalar multiwell potential, $U$, and the associated trajectory of \eqref{Langevin}.  As expected, the trajectory tends to persist near one of the two  minima.}
    \label{fig:langevin}
\end{figure}

Let us explain the rationale of the approach in a simpler, finite dimensional setting. Consider one-dimensional Langevin dynamics in multi-well potential:
\begin{equation}\label{Langevin}
dY_t = - U' (Y_t) dt + \sqrt{\alpha} d\beta_t
\end{equation}
where $Y_t \in \R, U:\R \rightarrow \R$ is a confining potential, $\alpha>0$, and $\beta_t$ is one dimensional  standard Brownian motion. Under very general  assumptions  on $U$, as soon as $\alpha>0$,  the SDE \eqref{Langevin} has only one invariant measure, the Maxwellian $e^{-2U/\alpha}$ (after choosing  $U$ to include the correct normalization constant, which can always be done). 
In contrast, when $\alpha=0$ this equation is a simple ODE with as many stable stationary states   as the number of minima (wells) of the potential $U$ (assuming such minima are isolated). Dirac measures centered around such stationary states can be viewed as the (multiple) stable invariant measures of the ODE.  For a double well potential, one has three stationary states, two stable and one unstable. 
When $U$ has many wells, a naive way to count them is to simulate \eqref{Langevin} for long enough and count the number of modes of the measure $e^{-2U/\alpha}$, as this measure will concentrate around the minima of $U$.  This is shown in Figure \ref{fig:langevin}.

This approach comes with limitations. First and foremost, one never knows if the whole state space has been adequately sampled; hence there are no guarantees that all the minima of $U$ have been found. This problem is unavoidable, whether in finite or infinite dimensions, and irrespective of the method used to sample. The potential landscape might also be more complex with non-isolated minima, elusive saddle points, etc. All these problems persist in the infinite dimensional setting. This is why we propose this method  to be used in conjunction with other methods, as an aid for exploration.  We also emphasize, again, that the method we propose will only help count the number of {\em stable} solutions of \eqref{eqn: stationary PDE}, in the sense that they are stable in \eqref{eqn: time-dependent PDE}. It will be far less informative about unstable solutions.  

Analogously with the finite dimensional case, one might see the idea of using the SPDE \eqref{eqn: additive SPDE intro} as a landscape exploration method, especially if the deterministic part of the equation has some gradient structure. Some comments on this - and words of caution - can be found in Note \ref{note: Maxwellian}. 

While infinite dimensional Markov Chain Monte Carlo (MCMC) sampling algorithms are now well-developed for many problems of interest,  one might not know the invariant measure of \eqref{eqn: additive SPDE intro} explicitly, making the use of MCMC algorithms impractical. The option - which we take in the example of this paper - is to simulate the SPDE for sufficiently long time.
Despite all the caveats that we have highlighted, this approach is {simple to understand and general, in the sense that, since it relies on a concept that does not hinge on the specific structure of the equation at hand,  it can be also applied to equations different from our test case. More comments in this in Section \ref{sec: background}.}  In the next section we provide an illustration of its effectiveness in the context of the search of stationary solutions of   porous media equations.

\section{Background on the McKean-Vlasov PDE and associated SPDE}\label{sec: background}
As test case we consider  the problem of finding the number of stationary solutions of the McKean-Vlasov (MKV) PDE,
\begin{align}\label{PDEsimpleparticlesystem}
\pa_t \rho_t(x)  = \pa_x \left[\left(V'(x) + (F'\ast \rho_t)(x)\right) \rho_t(x) + \sigma \pa_x\rho_t(x)\right] \,
\end{align}
for the unknown $\rho= \rho(t,x):\R_+ \times \T \rightarrow \R$, where $\T$ denotes the one dimensional torus given by the interval $[0, 2\pi]$  with periodic boundary conditions. \footnote{For time-dependent quantities we use indistinguishably the notation $X_t$ and $X(t)$, so throughout $\rho_t(x)=\rho(t,x)$. }
That is, we want to find the stable solutions of the problem \eqref{eqn: stationary PDE} when the operator $\cL$ is the non-linear differential operator given  by the right hand side of \eqref{PDEsimpleparticlesystem}, which we henceforth refer to as the {\em McKean-Vlasov differential operator}. 
In the above 
$V, F: \T \rightarrow \R $ are given functions, commonly referred to as the environmental and inter-particle potential, respectively,  $'$ denotes derivative with respect to the argument of the function and $\ast$ denotes convolution. Let us briefly recap what is known about the evolution \eqref{PDEsimpleparticlesystem}.   

The PDE \eqref{PDEsimpleparticlesystem} is a popular and flexible model in a range of applications and  has been studied for decades. It models the evolution of the density of a system of interacting particles. More precisely it is associated to the  following interacting particle system (IPS)
\be\label{initialPS}
dX_t\iin= - V'(X_t\iin) - \frac{1}{N}\sum_{j=1}^N F'(X_t\iin-X_t\jn)  dt+ \sqrt{2\sigma} d\beta_t^i, 
\ee
where, for every $i=1\dots N$, $X_t\iin \in \T$ is the position of the $i$-th particle at time $t\geq 0$. Using mean field techniques \cite{talay1996probabilistic}  one can show that the empirical measure $\mu^N_t:= \frac1N\sum_{i=1}^N \delta_{X_t\iin}$ of the system converges, as $N$ tends to infinity, to a deterministic probability measure. Such a measure has a density, which is precisely the solution $\rho$ of the PDE \eqref{PDEsimpleparticlesystem}. The PDE \eqref{PDEsimpleparticlesystem} can also be viewed as the Fokker-Planck equation of the following non-linear SDE
\begin{align}
dX_t & = - \left[V'(X_t) + \int_{\T} F'(X_t-y)\rho_t(dy)
\right] dt+ \sqrt{2\sigma} d\beta_t \,,  \quad X_t \in \T   \,.\label{nonlSDE}
\end{align}
In the above $\rho_t$ is the density of the law of $X_t$ -  so that the SDE is non-linear in the sense of McKean - and,  by It\^o's formula,  one can see that $\rho_t$ solves  \eqref{PDEsimpleparticlesystem}. Hence the  stationary solutions of \eqref{PDEsimpleparticlesystem} are the invariant measures of \eqref{nonlSDE}. 
A key fact to observe is that the evolution equation \eqref{PDEsimpleparticlesystem} preserves mass and positivity; that is, $\int_{\T} \rho_t(x) dx $ is a constant, independent of time, and $\rho_0\geq 0 $ implies $\rho_t\geq0$. This is coherent with the probabilistic interpretation of the PDE. 

The problem of determining the set of stationary solutions of \eqref{PDEsimpleparticlesystem} is far from trivial, even when the spatial variable is one-dimensional.  Following the seminal paper \cite{dawson1983critical}, this problem has been tackled in a variety of settings, \cite{bertini2014synchronization, pavliotis, constantin2004note, constantin2004remarks, herrmann2010non}. In particular, the number of stationary solutions of \eqref{PDEsimpleparticlesystem}  depends on the choice of potentials $V$ and $F$ and on the strength of the noise.  
Here, to fix ideas, we consider the case where $V$ and $F$ are given by
\begin{equation}\label{eqn:torusVF}
V(x) = \cos (2x), \quad F(x)= -\cos(x) \,.
\end{equation}
The relevant properties of this choice is that  $V$ is non-convex (more specifically, double well and symmetric) and $F$ introduces an attractive force between particles. One can show that with this choice of $V$ and $F$, there exists a critical value of the noise, $\sigma_c>0$,  such that  if $\sigma>\sigma_c$ equation \eqref{PDEsimpleparticlesystem} admits exactly one stationary solution, $\rho_0$, which is bimodal and  equally concentrated around the bottom of  two wells. The two wells have minima at $x_1=\pi/2$ and $x_2=3\pi/2$.  For $1/2<\sigma<\sigma_c$ and for $\sigma$ small enough there are exactly three stationary solutions of \eqref{PDEsimpleparticlesystem} - the conjecture being that for any $0<\sigma\leq \sigma_c$ there should be three stationary solutions. One is $\rho_0$ itself; the other two, $\rho_{\pm}$, are still  bimodal with modes around $x_1$ and $x_2$ but $\rho_+$ concentrates most of the mass in  $x_1$ while $\rho_-$ concentrates most of its mass around $x_2$.   The reader may find intuition for this property in, amongst other references, \cite{angeli2023well}. 
Here we briefly recall (as it will be relevant to Section \ref{sec:simulations} and for later comments) that this result is obtained by first observing that, irrespective of the choice of $V$ and $F$, any stationary solution $\rho_{\infty}$ of \eqref{PDEsimpleparticlesystem} solves the following fixed point problem
\begin{equation}\label{fixed point eqn} 
 \rho_{\infty}(x) = \frac{1}{Z_{\sigma}}
 e^{-\frac 1\sigma \left[  V(x) + \int \limits_0^x \left( F^{'} \ast \rho_{\infty} \right)(y) dy \right]} \,.
\end{equation}
Once  $V$ and $F$ are chosen as in \eqref{eqn:torusVF}, the above, a priori, infinite dimensional fixed point problem, can be reduced to a finite-dimensional, two-parameter fixed point problem. One can indeed observe (see \cite[Section 2]{angeli2023well}) that in this case $\rho_{\infty}$ solves \eqref{fixed point eqn} if and only if it has the following structure
\begin{equation}\label{eqn: stat solution general eqn}
    \rho_{\infty}(x)=\frac{1}{Z_{\sigma}}
 e^{-\frac 1\sigma \left[  \cos\,2x -m_1\sin\,x-m_2\cos\,x \right]},
\end{equation}
where $m_1,m_2 \in \R$ are fixed point of the following equations
\begin{equation}\label{eqn:two parameter fixed point problem}
    \int \limits_{\T}\rho_{\infty}(x)\,\sin\,x\,dx = m_1,\quad \int \limits_{\T}\rho_{\infty}(x)\,\cos\,x\,dx=m_2 \, ,
\end{equation}
respectively. To recap, when $\cL$ is the McKean-Vlasov differential operator, solving \eqref{eqn: stationary PDE} corresponds to solving the fixed point problem \eqref{fixed point eqn}. The results we have summarised above on the number of stationary states of the McKean-Vlasov PDE are obtained by using this strategy. More comments on this in Note \ref{note: Maxwellian}. 

We note in passing that the fact that the PDE \eqref{PDEsimpleparticlesystem} can admit several stationary solutions is in contrast with the behaviour of the  particle system \eqref{initialPS}, which, irrespective of the strength of the noise $\sigma>0$, always admits just one invariant measure (although clearly this invariant measure is made of multiple metastable states). So the particle system is not necessarily helpful to our goals (see Note \ref{note: Maxwellian}). 

While \eqref{PDEsimpleparticlesystem} has many stationary states, 
it is natural to imagine that, after adding a strong-enough noise to \eqref{PDEsimpleparticlesystem}, we will obtain an SPDE with a unique invariant measure (this invariant measure is a measure on the space of measures). Technical details aside, this turns out to be the case. Indeed, consider the SPDE
\begin{align}\label{eqn:additiveSPDE}
\pa_t u_t(x)  = \pa_x \left[\left(V'(x) + (F'\ast u_t)(x)\right) u_t(x) + \sigma \pa_x u_t(x)\right] + Q^{1/2} \pa_t W_t
\end{align}
where $W_t$ is cylindrical Wiener noise and $Q$ is a trace class operator, so that 
\be\label{cylnoise}
Q^{1/2}\partial_t W_t = \sum_k\lambda_k e_k(x) dw_t^k, 
\ee
and the $w_t^k$'s are independent one dimensional  standard Brownian motions, $\{e_k\}_k$ is a basis of $L^2(\T;\R)$ (in which $Q$ is assumed to diagonalise) and  $\{\lambda_k^2\}_{k\in \N}$ is the sequence of non-negative eigenvalues of $Q$. 
We will assume that the eigenvalues of $Q$ are of the form
\begin{equation}\label{eqn: eigenvalues of  Q}
  \lambda_0=0 \, \, \mbox{ and } \, \,   \lambda_k^2 = \frac{\gamma^2}{k^{2s}}e_k, \mbox{ for } k\geq 1
\end{equation}
for some $\gamma>0$ and $s>1/2$.
Well-posedness and ergodicity of the SPDE \eqref{eqn:additiveSPDE} have been studied in \cite{angeli2023well}. Using the results of \cite{angeli2023well} one can easily see that with the above choice of $Q$, the SPDE \eqref{eqn:additiveSPDE} is well-posed in $L^2(\T;\R)$ provided $s > \frac{1}{2}$, while it is strong Feller and irreducible, respectively,  for $\frac{1}{2}<s<1$ and $s>\frac{1}{2}$, respectively. This implies that for $1/2<s<1$ the SPDE admits at most one invariant measure.  Whether such an invariant measure exists or not is still an open problem, but we believe it to be the case. 
 We proceed by assuming \eqref{eqn:additiveSPDE} admits a unique invariant measure, at least under \eqref{eqn: eigenvalues of  Q}.  Our simulations in Section \ref{sec:simulations} will, as a by-product, support this conjecture. 

To our purposes a more pressing issue is the fact that the  SPDE \eqref{eqn:additiveSPDE} does not preserve mass in general and it does not preserve positivity either. This means that sampling from the invariant measure of \eqref{eqn:additiveSPDE} may also result in exploring parts of state space in which we are not interested, as we only want stationary solutions of \eqref{PDEsimpleparticlesystem} which have mass one and are non-negative.  The choice $\lambda_0=0$ in \eqref{eqn: eigenvalues of  Q} has been made to remedy this issue  and  guarantee mass conservation - by formally integrating the equation in space and using the fact that $\int e_k(x) \, dx=0$ for every $k\geq 1$ one can see that, if $\lambda_0=0$, then the quantity $\int u_t(x)\,dx$ is constant in time. As for the positivity, a straightforward calculation shows that stationary solutions of \eqref{PDEsimpleparticlesystem} with mass one are always non-negative.\footnote{{The proof of the non-negativity of stationary solutions of \eqref{PDEsimpleparticlesystem} is straightforward. One first takes integrates w.r.t. $x$ the second order ODE  $\pa_x \left[\left(V'(x) + (F'\ast \rho_{\infty})(x)\right) \rho_{\infty}(x) + \sigma \pa_x\rho_{\infty}(x)\right]$\,=\,0. After that, one solves via variations of constants the resulting first order ODE. Then imposes the periodic boundary conditions (PBC) on the solution and the fact that it is a probability measure and finally the expression \eqref{eqn: stat solution general eqn} is obtained. It is clear that functions with the expression such as \eqref{eqn: stat solution general eqn} are non-negative.}} Therefore, once condition \eqref{eqn: eigenvalues of  Q} is enforced, we expect the solution of \eqref{eqn:additiveSPDE} to concentrate only around the desired (stable) stationary solutions of \eqref{PDEsimpleparticlesystem}, i.e. those which have mass  mass one  and are hence positive.

To the best of our knowledge, when $0<\sigma<\sigma_c$, there are no results on  the global stability of the three stationary states  of \eqref{PDEsimpleparticlesystem} (with $V$ and $F$ as in \eqref{eqn:torusVF})  - though some partial results in this direction have been obtained in \cite{tugaut2013convergence} and later in \cite{bashiri2020long}.  When $0<\sigma<\sigma_c$ one expects that two of the equilibria ($\rho_{\pm}$) will be stable and one unstable ($\rho_0$). In Appendix \ref{appendix:stability} we numerically show that this is indeed the case.   Hence, when sampling from the invariant measure of \eqref{eqn:additiveSPDE}, we expect to see bimodal behavior in our measure. This is confirmed by the simulations of Section \ref{sec:simulations}; more detail is provided in that section.  In general the linearly unstable solutions may be akin to saddle points; they will be challenging to see directly in the simulation of \eqref{eqn:additiveSPDE}, but may be critical to understanding transition rates amongst metastable solutions.

An explicit form for the invariant measure of \eqref{eqn:additiveSPDE} is not known. This precludes the use of many MCMC sampling strategies as they often depend on evaluating the (unnormalized) invariant density in a Metropolis ratio. {Moreover, since we are looking to sample from a measure supported on an infinite dimensional space, it would be advisable to use an MCMC method which is well posed in infinite dimension \cite{cotter2013mcmc, hairer2005analysis, hairer2007analysis}. But, again, we do not know of any MCMC method which is well posed in infinite dimension and that can be implemented in absence of information on the analytical form of the target measure \cite{cotter2013mcmc, beskos2008mcmc, dobson2020reversible}.  }

We therefore take the approach of simulating the SPDE \eqref{eqn:additiveSPDE} for sufficiently long.   As the simulation of this SPDE needs care, we present in Section \ref{sec:numerical discretization} a numerical scheme that is used for its study.  We reserve for future work the proof of the convergence of this numerical scheme. 

\begin{note}\label{note: Maxwellian}
\textup{Some further comments on the SPDE \eqref{eqn:additiveSPDE} and its simulation.}
\begin{itemize}
    \item \textup{The SPDE \eqref{eqn:additiveSPDE}  can be obtained as a limit of an appropriate (weighted) interacting particle system. This has been proved in   \cite{angeli2024approximation}. Such a particle system could be used as a way of simulating the SPDE.  However, similarly to what happens for the particle system \eqref{initialPS} and its limiting PDE \eqref{PDEsimpleparticlesystem}, it may be the case that the ergodic properties of the SPDE \eqref{eqn:additiveSPDE} are different from those of the associated particle system. Indeed, the work \cite{angeli2024approximation} proves that the particle system converges to the SPDE over finite time horizons, but there is a priori no reason to believe that convergence will hold uniformly in time, i.e. over infinite time horizons.  This is why we opted to simulate the SPDE with a numerical discretization scheme.}
    
    \item \textup{As we mentioned at the end of Section \ref{subsec: approach intro},  the  deterministic part of the evolution \eqref{eqn:additiveSPDE} is in gradient form (see \cite{pavliotis}). So one would expect, in analogy to the finite dimensional setting discussed around equation \eqref{Langevin}, that the approach we take here could be interpreted as a landscape exploration method. Heuristically this is perhaps the case. However, there is a subtlety to this point. As we mentioned, there are no rigorous results on the existence of the invariant measure of \eqref{eqn:additiveSPDE}. In particular, we do not know if such a measure is of Maxwellian type, and this may be challenging, in view of recent observations \cite{delarue2022rearranged}. Interpreting the method we consider here as an infinite dimensional  landscape exploration method hinges on making progress on this aspect of the theory. Nonetheless,   whether the invariant measure is of Maxwellian type or not, one still expects that such a measure will concentrate around the stable stationary solutions of the McKean-Vlasov PDE (as confirmed by our simulations), and this is all we truly need. } 
    
    \item{\textup{Still, on the matter of potentially using particle systems to solve the problem of interest, it should be noted that the particle system \eqref{initialPS} has, for each fixed $N$, a unique invariant measure. Such an invariant measure is metastable and the number of metastable states of this measure coincides with the number of stable stationary states of \eqref{PDEsimpleparticlesystem} \cite{berglund2019introduction} - and it is expected to be the same as the number of metastable states of the invariant measure of the SPDE \eqref{eqn:additiveSPDE}. However the exit time from a potential well for the particle system \eqref{initialPS} grows exponentially with $N$ \cite{berglund2019introduction, herrmann2008large}, making this approach computationally impractical (if not infeasible), especially in the presence of many wells. In contrast, using the metastability of the SPDE \eqref{eqn:additiveSPDE} is computationally much faster.}}
    
    \item \textup{Reducing the infinite dimensional fixed point problem \eqref{fixed point eqn} to a finite dimensional one is only practical when $V$ and $F$ are sufficiently simple. If $V$ and $F$ allow for this drastic reduction in dimensionality (with our choices of potentials this is due to the fact that such potentials are each given by a single Fourier mode),  then it is usually possible to solve, numerically, the resulting finite dimensional fixed point problem. For us, this is problem \eqref{eqn:two parameter fixed point problem}. However, even in these simple cases, the analytical study of the finite dimensional fixed point problem is challenging, and it tends to offer only partial results (as for the case we have detailed above).  See  \cite{angeli2023well}, references therein  and the list of references at the start of this section.   \\
    Even slightly complicating  potentials $V$ and $F$ makes the reduction of the infinite dimensional fixed point problem to a finite dimensional one substantially more difficult. Should it succeed, the numerical  and  analytical study, respectively,  of the resulting  finite dimensional problem  become tricky and  prohibitive to pursue. In Section \ref{sec:simulations} we choose $V$ and $F$ as in \eqref{eqn:torusVF}; this provides a good test case, as the number of stable stationary states is  known in this case (modulo the comments we made above).  However, applying our method is not needed in this relatively simple case. \\
    We expect  our method will be more useful in more complex situations - including complex choices of $V$ and $F$ or, more so, when the PDE part of the evolution is not of McKean-Vlasov type. Indeed, if the PDE part of the evolution is not of McKean-Vlasov type then  it is not necessarily the case that  the problem \eqref{eqn: stationary PDE} can be reduced to some form of fixed point equation.  In such instances, the method we have present  allows for greater flexibility. An interesting avenue could be the one of using the method we present here for different non-linear PDEs (using the associated SPDEs), such as Allen-Cahn type equations, Navier stokes, and equations appearing in non-linear elasticity, such as those in \cite{wang2019order}. This will however require new theoretical efforts in the ergodic theory for SPDEs. }
\end{itemize} 
\end{note}

\section{Numerical Methods}
\label{sec:numerical discretization}
Here, we review our computational methods for exploring \eqref{eqn:additiveSPDE}.

\subsection{Discretization of the SPDE}
To generate sample paths of \eqref{eqn:additiveSPDE}, we employ a spectral Galerkin spatial discretization together with Euler-Maruyama time stepping.  We refer the reader to, for instance, \cite{lord_introduction_2014}. Denoting the total number of resolved Fourier modes $J$, and using complex trigonometric functions, we will achieve spectral accuracy with $\mathrm{O}(J\ln J)$  scaling through the Fast Fourier Transform (FFT).  Euler-Maruyama time stepping for this problem (having only additive noise) will yield a strong error of $\mathrm{O}(\Delta t)$ in the sample paths.  This strategy is entirely standard, and though we do not provide a full numerical analysis of it in the present work, our numerical tests are consistent with what would be expected from this method. We defer such numerical tests to Section \ref{sec:simulations} \ref{subsec:verification of algorithm}. We expect that a complete analysis of the numerical scheme we present here should be accessible and we hope to pursue this  in future work. 

While our example computations will be for the case that $F$ and $V$ are given by \eqref{eqn:torusVF}, these methods are expected to work for all sufficiently regular $V$ and $F$. 

\subsubsection{Space Discretization: Spectral Galerkin Method.}
The fundamental assumption is that the covariance, $Q$, of the noise process in \eqref{eqn:additiveSPDE} diagonalizes against the the complex trigonometric functions on the $\mathbb{T}=[0,2\pi)$  
\begin{equation}
\label{e:trig}
e_k = \frac{1}{\sqrt{2\pi}}e^{i k x}, \quad Q e_k = \lambda_k^2  e_k, \quad k \in \mathbb{Z}.
\end{equation}
This us allows to express the noise as \eqref{cylnoise}.  We now represent the solution, $u$, of \eqref{eqn:additiveSPDE} by
\begin{equation}
    u(x,t) = \sum_{k} \hat{u}_k(t)e_k.
\end{equation}
Substituting in and projecting onto $e_k$, \eqref{eqn:additiveSPDE} becomes
\begin{equation}
    \hat{u}_k' = ik \left(  \mathcal{F}[V'\cdot u](k)  + \mathcal{F}[\mathcal{F}^{-1}[\hat{F'}\cdot\hat{u}]u](k)\right) - \sigma k^2 \hat{u}_k + {\lambda_k} w'_k
\end{equation}
Symbols $\mathcal{F}$ and $\mathcal{F}^{-1}$ denote the Fourier transform and its inverse; in software, these are carried out by the FFT algorithm.  The above notation is to emphasize that $V'\cdot u$ is computed in real space, while $F'\ast u$  is first computed in Fourier space, and then $(F'\ast u)\cdot u$ is computed in real space.


Using an integrating factor on the linear term, we can write the mild solution (in Fourier coordinates) as
\begin{equation}
\label{e:specdisc1}
\begin{split}
    \hat{u}_k(t) &= e^{-\sigma k^2 t}\hat{u}_k(0) \\
    &+ \int_0^t i k e^{-\sigma k^2 (t-s)}\left( \mathcal{F}[V' u](k) 
    + \mathcal{F}[\mathcal{F}^{-1}[\hat{F'}\cdot\hat{u}]u](k)\right)(s)ds + {\lambda_k}\int_0^t e^{-\sigma k^2 (t-s)}dw_k
\end{split}
\end{equation}
Thus far, we have only formally manipulated \eqref{eqn:additiveSPDE}; no approximation has  been made.  We now truncate to $J$ (assumed to be even) modes, so that
\begin{equation}
\label{e:uJ}
    u^{(J)}(t) = \sum_{k=-J/2+1}^{J/2} \hat{u}_{k}^{(J)}e_k(x)
\end{equation}
and let $P^{(J)}$  be the projection onto the $J$ lowest modes.  Substituting \eqref{e:uJ} into \eqref{e:specdisc1} and applying this projection we obtain our finite dimensional system:
\begin{equation}
\label{e:specdisc2} 
\begin{split}
    \hat{u}_{k}^{(J)}(t) &= e^{-\sigma k^2 t}\hat{u}_{k}^{(J)}(0) \\
    &+ \int_0^t i k e^{-\sigma k^2 (t-s)}P^{(J)}\left( \mathcal{F}[V' {u}^{(J)}](k) 
    + \mathcal{F}[\mathcal{F}^{-1}[\hat{F'}\cdot\hat{u}^{(J)}]{u}^{(J)}](k)\right)(s)ds \\
    &+ {\lambda_k}\int_0^t e^{-\sigma k^2 (t-s)}dw_k
\end{split}
\end{equation}
As is sufficient for problems with a quadratic nonlinearity, the projection operator can be implemented in code by:
\begin{itemize}
    \item Simulating with $2J$ modes, if we wish to resolve $J$ of them;
    \item Zeroing out the $J$ highest frequency modes.
\end{itemize}
This corresponds to anti-aliasing, \cite{boyd_chebyshev_2001}. 

\subsubsection{Time Discretization: Euler-Maruyama}
The fully discretized problem is obtained via an Euler-Maruyama temporal discretization of \eqref{e:specdisc2}. Hence, on top of the spatial approximation, we make the temporal approximation
\begin{equation}
\hat{u}^{(J, \Delta t)}_{k,n} \approx \hat{u}_k^{(J)}(t_n),
\end{equation}
with $t_n = n \Delta t$. After a bit of stochastic calculus on the noise term, the update is:
\begin{equation}
\label{e:fullydisc1}
\begin{split}
     \hat{u}^{(J,\Dt)}_{k,n+1} &=  e^{-\sigma k^2 \Delta t}\hat{u}^{(J,\Dt)}_{k,n} \\
     &+ ik\frac{1-e^{-\sigma k^2 \Delta t}}{\sigma k^2}P^{(J)}\left(\mathcal{F}[V' {u}^{(J,\Dt)}] 
    + \mathcal{F}[\mathcal{F}^{-1}[\hat{F'}\cdot\hat{u}^{(J,\Dt)}]{u^{(J,\Dt)}}]\right)_{k,n}\\
    &  + \lambda_k\sqrt{\frac{1 - e^{-2\sigma k^2\Delta t}}{2\sigma k^2}} \xi_{k, n+1}, \quad \xi_{k, n+1}\overset{\mathrm{i.i.d.}}{\sim} N(0,1).
\end{split}
\end{equation}


\section{Simulations}\label{sec:simulations}

\begin{figure}[!ht]
\centering\includegraphics[width=3.5in]{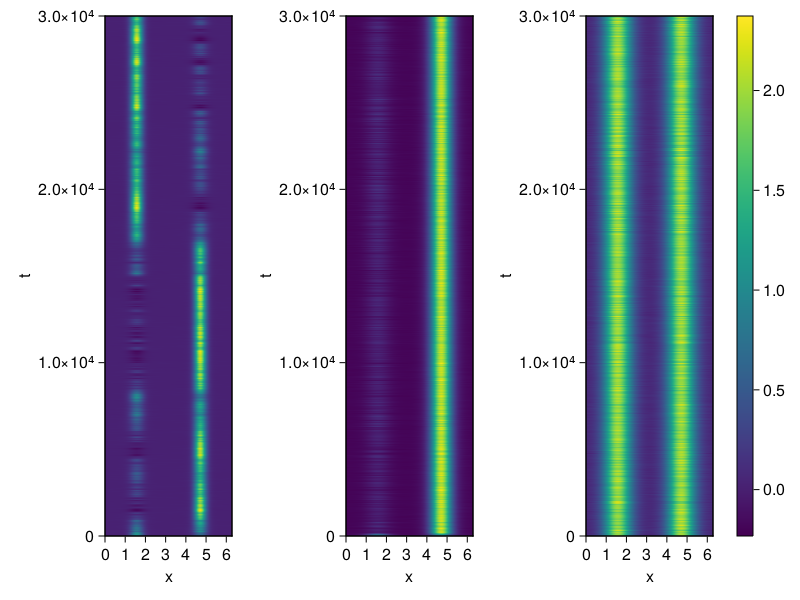}
\caption{Heat map of the solution $u_t$ to the SPDE \eqref{eqn:additiveSPDE} with $V$ and $F$ as in \eqref{eqn:torusVF}. From left to right, these plots have been obtained with $\sigma=0.2,0.6$ and $1$, respectively; we recall that $\sigma=1$ is above the critical threshold.   The other parameter values held fixed given by: $t_{\text{max}}=3\cdot10^4,\Delta t=10^{-2}, s=0.75, \gamma = 10^{-2}$ and initial condition $u_0(x)=\frac{1}{\pi}\sin^2\,x$. On the vertical axis of each graph there is time $t$ while on the horizontal axis lies the spatial variable $x$. The colour bar is  common across the three plots. These plots have been obtained by employing CairoMakie \cite{DanischKrumbiegel2021}.}
\label{fig_heat map double well}
\end{figure}

\begin{figure}[!ht]

\begin{subfigure}
    \centering
    \includegraphics[width=0.7\textwidth]{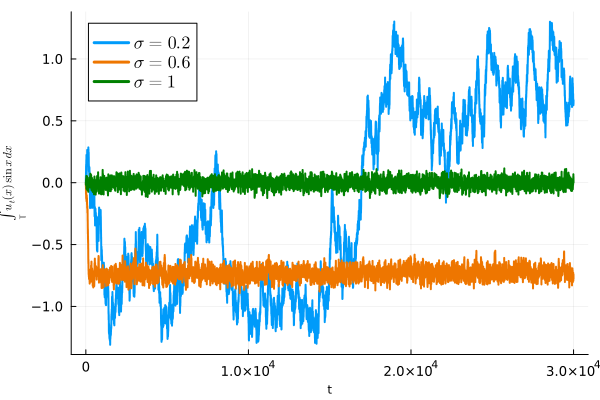}
    \caption{Plot of the function $t \to I_1(t)=\int \limits_{\T} u_t(x)\sin \,x\,dx$ for the solution $u_t$ to \eqref{eqn:additiveSPDE} with $\sigma=0.2,0.6$ and $1$.  These simulations are run with the same parameters as in Figure \ref{fig_heat map double well}.}
    \label{functional:I1}
\end{subfigure}
\begin{subfigure}
    \centering
    \includegraphics[width=0.7\textwidth]{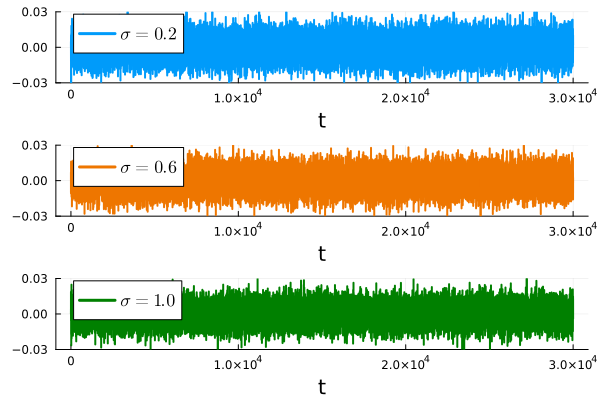}
    \caption{Plot of the function $t \to I_2(t)=\int \limits_{\T} u_t(x)\cos \,x\,dx$ for the solution $u_t$ to \eqref{eqn:additiveSPDE} with $\sigma=0.2,0.6$ and $1$.  These simulations are run with the same parameters as in Figure \ref{fig_heat map double well}.}
    \label{functional:I2}
\end{subfigure}
\end{figure}


\begin{figure}
\centering
\includegraphics[width=3in]{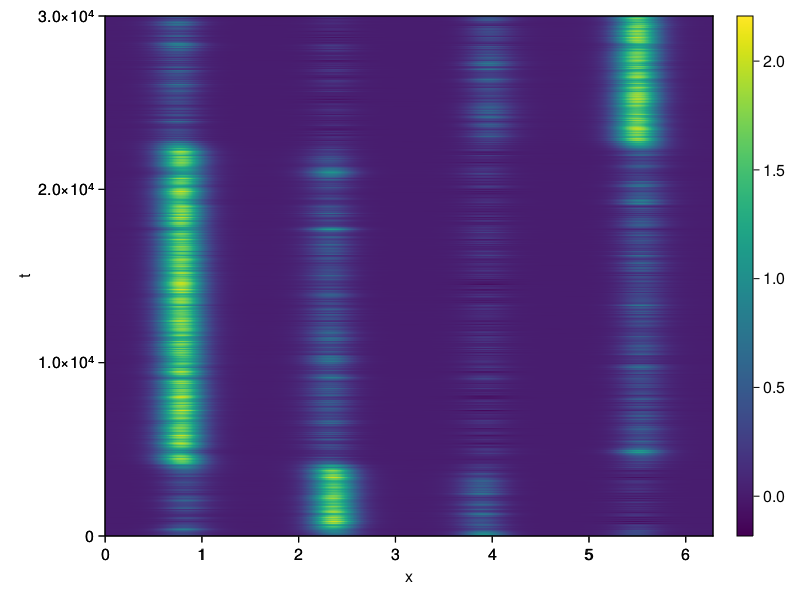}
\caption{Heat map of the solution $u_t$ to the SPDE \eqref{eqn:additiveSPDE} with $V(x)=\cos \,4x$, and $F$ as in \eqref{eqn:torusVF}. This plot has been obtained with the following choice of the parameter values: $\sigma=0.4, t_{\text{max}}=3\cdot10^4,\Delta t=10^{-2}, s=0.75, \gamma = 10^{-2}$ and initial condition $u_0(x)=\frac{1}{\pi}\sin^2\,x$. On the vertical axis  there is time $t$ while on the horizontal axis lies the spatial variable $x$.  The plots have been obtained by employing CairoMakie \cite{DanischKrumbiegel2021}.}
\label{fig_heat map four well}
\end{figure}

\subsection{Computed Stationary Solutions}

We now use the numerical discretization method presented in Section \ref{sec:numerical discretization} to simulate the dynamics \eqref{eqn:additiveSPDE} for long times. When the stationary solutions of \eqref{PDEsimpleparticlesystem} are isolated, we expect the process $\{u_t\}_{t\geq 0}$ solution of \eqref{eqn:additiveSPDE} to exhibit metastable behaviour, hopping between different stable stationary solution of \eqref{PDEsimpleparticlesystem}. To observe this fact numerically, we show  a heat plot of the solution $u_t$, see Figure \ref{fig_heat map double well}; moreover, inspired by   \eqref{eqn:two parameter fixed point problem}, we consider the quantities 
\begin{equation}
    I_1(t):= \int_{\T} u_t(x) \sin x \, dx, \quad  I_2(t):= \int_{\T} u_t(x) \cos x \, dx \,.
\end{equation}
When $V, F$ are as in \eqref{eqn:torusVF} and $\sigma$ is either small enough or $1/2<\sigma<\sigma_c$, we know that \eqref{PDEsimpleparticlesystem} has two stable stationary states, $\rho_{\pm}$. Such stable states correspond to the values of $(m_1,m_2)$ given by the points $(\pm 1,0)$ (for detail see \cite{angeli2023well}).  We expect $I_1(t)$ to be metastable around the values of  $m_1=\pm 1$ and $m_2$ to fluctuate around the origin.  This is consistent with {Figures \ref{functional:I1} and \ref{functional:I2}, respectively}. When $\sigma>\sigma_c$, which is the case for the right most plot in Figure \ref{fig_heat map double well}, we know that \eqref{PDEsimpleparticlesystem} has a unique stationary solution, which we expect to be bimodal and equally concentrated around the two minima.  This is confirmed by the simulation.    

\textup{The value of $\gamma=10^{-2}$ used in these experiments was found empirically.  As we had some foreknowledge of the stationary solutions in the system, we found this value to be  sufficiently small value such that these solutions were metastable on the timescale of the simulation.  It is possible that we may miss solutions should this noise level be too large.}

Let us now make a further remark about Figure \ref{fig_heat map double well}. The left ($\sigma=0.2$) and right ($\sigma = 1$) panels  in that figure are consistent with our expectations. Since the value $\sigma=0.6$ is also below the critical threshold, we would have expected to see metastability also in the panel in the middle. That is, we would have expected the central panel to look a bit more like the one on the left. Upon running the simulation for longer, the expected metastability appears also for the value $\sigma=0.6$ (not shown in figures). Nonetheless  it is at first sight counter intuitive that a higher value of $\sigma$ should be slowing down  the hopping between metastable states. We do not have a complete explanation of this phenomenon yet, and we intend to carry out further study of this. However, we point out that this result is not  wholly unreasonable: in this infinite dimensional scenario the frequency of hopping between metastable states is governed by the parameter $\gamma$, which controls the intensity of the infinite dimensional noise,  not by $\sigma$.  Increasing $\gamma$ does indeed  result in faster hopping between metastable states (not shown in figures). We do believe that understanding why an increase of $\sigma$ slows down metastability is related to the role that $\sigma$ plays in the (weighted) interacting particle system  associated to \eqref{eqn:additiveSPDE} - see \cite{angeli2024approximation}. 

Finally, in Figure \ref{fig_heat map four well} we run the simulation for different potentials $V$ and $F$. In particular the potential $V$ we choose in Figure \ref{fig_heat map four well} has four wells. For these potentials no rigorous proofs are available, though the result we find is in line with heuristics.

\subsection{Verification of Algorithm}\label{subsec:verification of algorithm}
\begin{figure}
    \subfigure[Convergence in $J$ with $\Dt = 10^{-2}$]{\includegraphics[width=6cm]{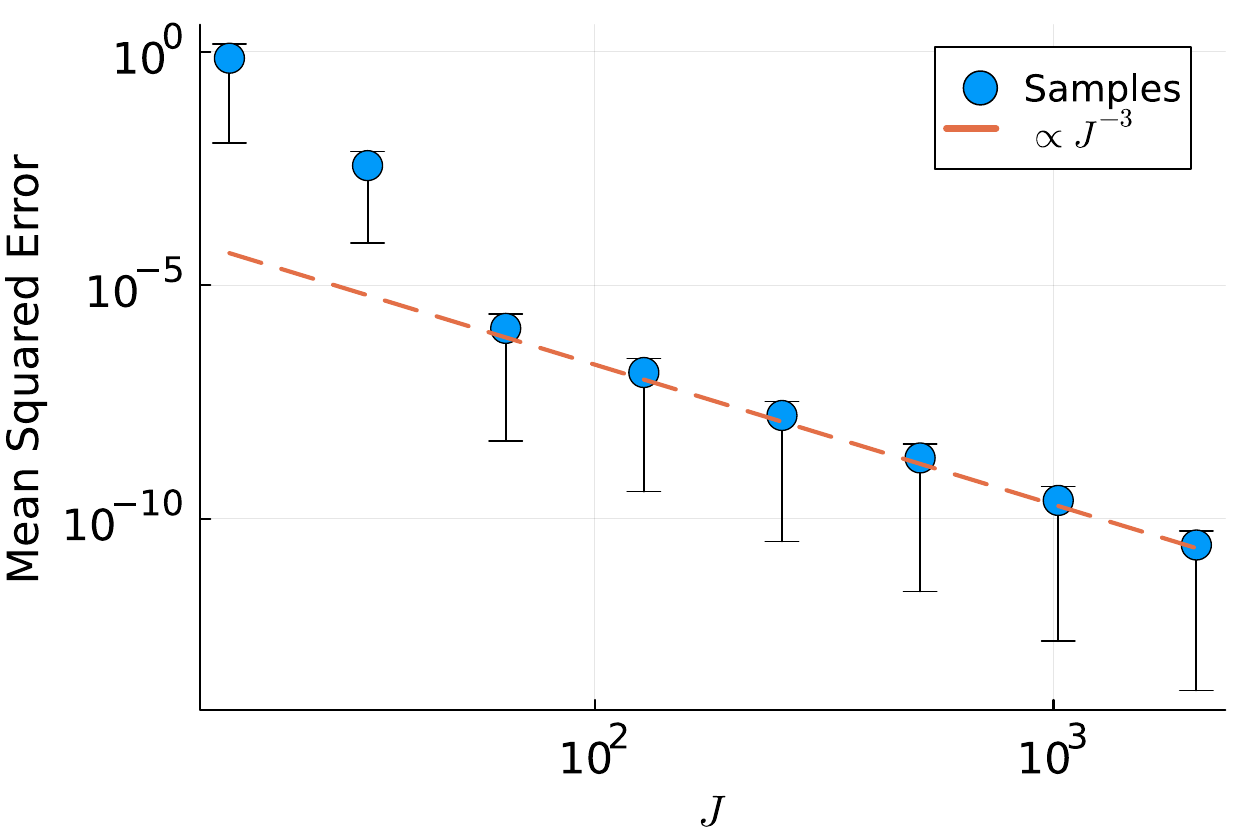}}
    \subfigure[Convergence in $\Dt$ with $J = 128$]{\includegraphics[width=6cm]{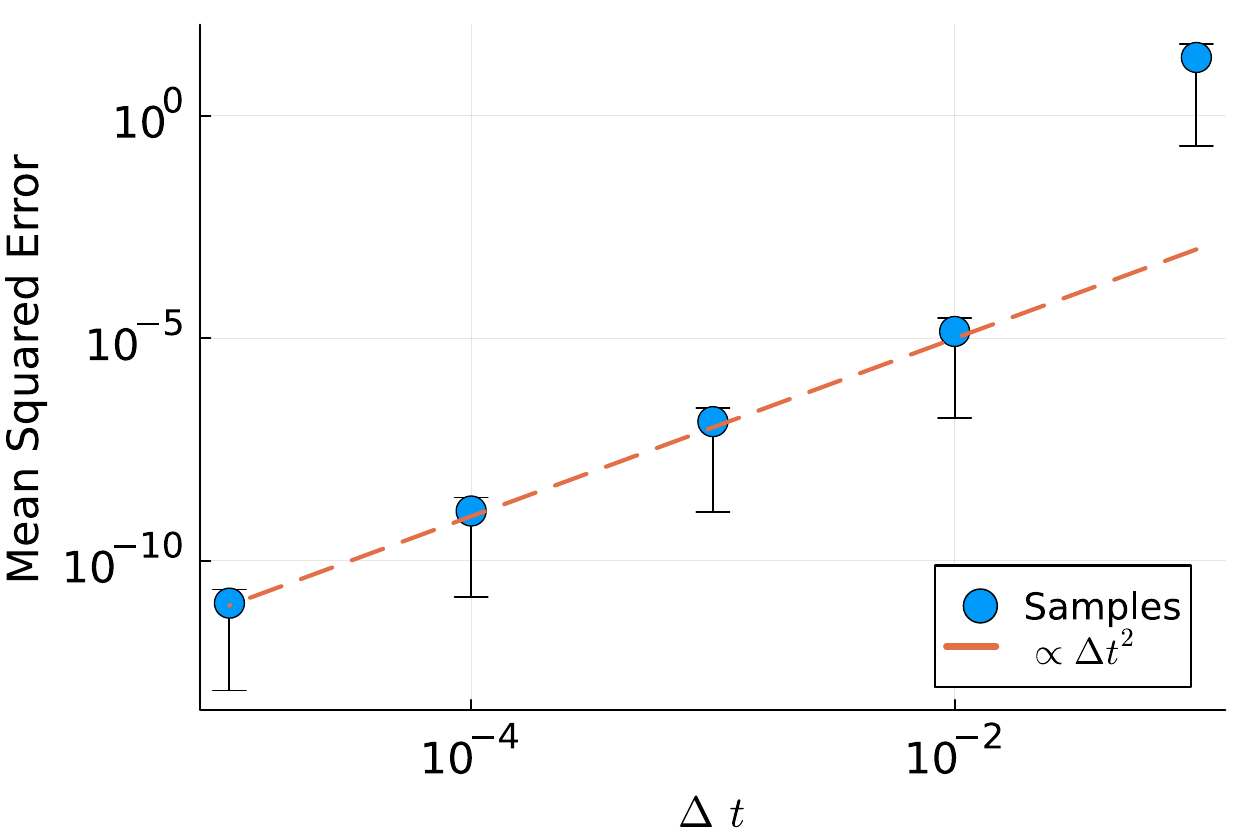}}

    \caption{Empirical verification of our algorithm's performance.  We see numerical convergence in both $J$ and $\Dt$ at the anticipated rates.  Computed with $\sigma=\gamma=0.1$, $s=1$, integrated until $t_{\max}=10$.  Error bars reflect 95\% bootstrap percentile confidence intervals generated from $10^4$ independent trials.}
    \label{fig:conv}
\end{figure}

While a full numerical analysis is beyond the scope of the present work, we did perform some basic convergence testing for the algorithm (numerical discretization) presented in Section \ref{sec:numerical discretization},  see  Figure \ref{fig:conv}.  \textup{These plots display the convergence of an estimate of the Mean Squared Error
\begin{equation}
    \E[\|u^{(J,\Delta t)}(t_{\max})- u(t_{\max})\|_{L^2}^2]\approx \frac{1}{n_{\rm trials}} \sum_{i=1}^{n_{\rm trials}}\|u^{(J,\Delta t, i)}(t_{\max}) -u^{\rm ref}(t_{\max})\|_{L^2}^2
\end{equation}
as we vary $J$ and $\Delta t$ independently, with $n_{\rm trials}= 10^4$ independent trials in each case.  The $L^2$ norms are computed in Fourier space.  For reference solutions, we used: $\Dt_{\rm ref} = 10^{-2}$ and $J_{\rm ref}=2^{12}$ to check $J$ convergence; $\Dt_{\rm ref } = 10^{-6}$ and $J_{\rm ref}=2^7$ to check $\Delta t$ convergence.}  The convergence rates of both the number of Fourier modes $J$ the and time step $\Dt$ is consistent with analogous discretizations of semilinear SPDE (i.e., stochastic Allen-Cahn), \cite{lord_introduction_2014}.

\appendix\section{Stability of stationary states of \eqref{PDEsimpleparticlesystem}}\label{appendix:stability}

To investigate the stability of the stationary solutions to the deterministic PDE  we linearise the McKean-Vlasov operator around the equilibria (see \cite[Subsection 3.2]{pavliotis}) and study its resulting spectrum. Namely, let $\rho_{\infty}$ be a stationary solution to \eqref{PDEsimpleparticlesystem} then the linearised McKean-Vlasov operator around the equilibrium position $\rho_{\infty}$ is given  by 
\begin{equation}\label{linearised}
    L\eta = \sigma \pa_{xx}\eta+ \pa_x \left [V^{'}\eta+(F^{'}*\eta) \rho_{\infty}+(F^{'}*\rho_{\infty})\eta\right ],
\end{equation}
where $\eta \in L^2(\T;\R)$ with $\int\eta(x)\,dx=0$.  Under the same kind of discretization in $x$ as used in the main body of the paper, we thus study the eigenvalue problem with matrix $\mathbf{L}$,
\begin{equation}
\label{disclinearised}
    \mathbf{L} =\sigma \mathbf{D}_x^2 + \mathbf{D}_x\left(\diag(V'(\mathbf{x})+\mathbf{F}'\boldsymbol{\rho}) + \diag(\boldsymbol{\rho})\mathbf{F}'\right)
\end{equation}
The spectral differentiation matrix, $\mathbf{D}_x$ is dense Toeplitz,\cite{trefethen_spectral_1997}.  $\mathbf{F}'$ corresponds to the discretization of $F' \ast$, and $\boldsymbol{\rho}$ is the discrete approximation of $\rho_\infty$.  For the values of $J$ we consider, a standard eigenvsolver suffices.  

\begin{figure}
\subfigure[$\sigma=0.2$]{\includegraphics[width=4cm]{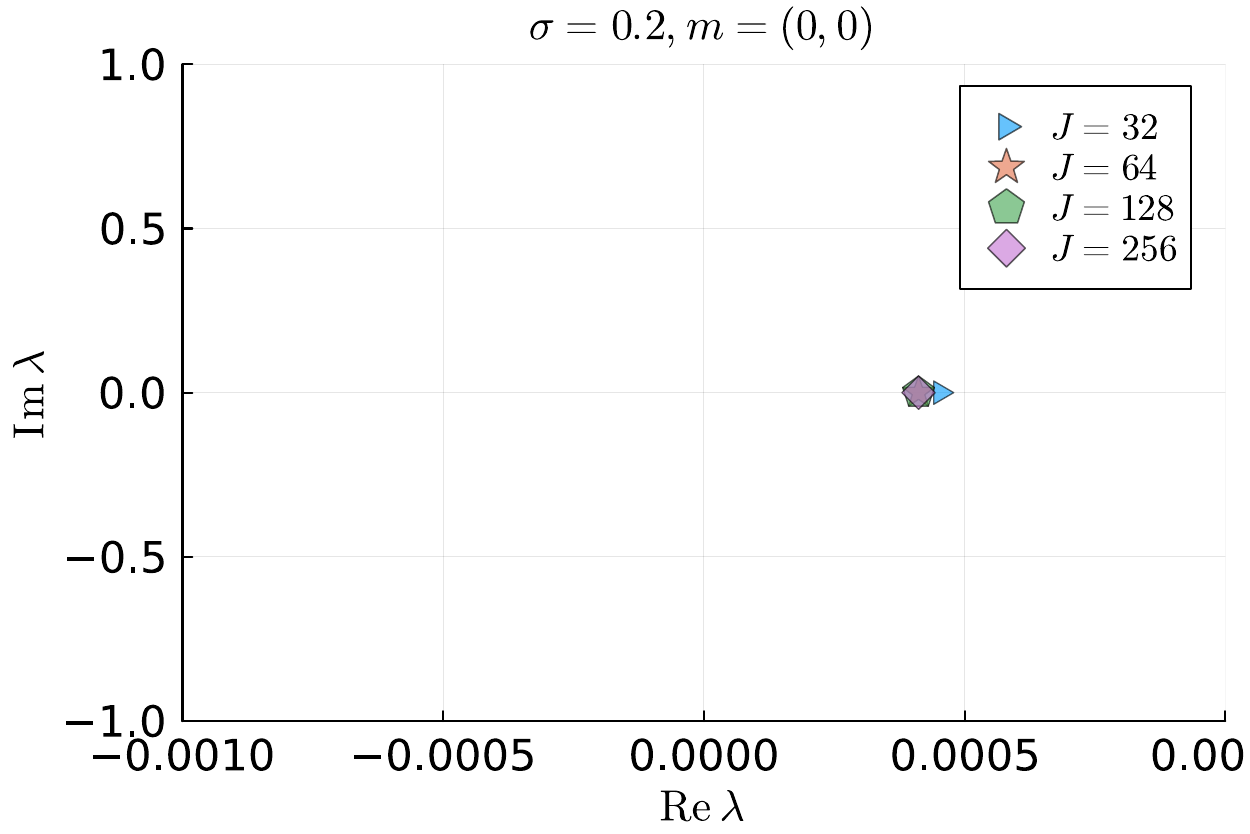}}
\subfigure[$\sigma=0.6$]{\includegraphics[width=4cm]{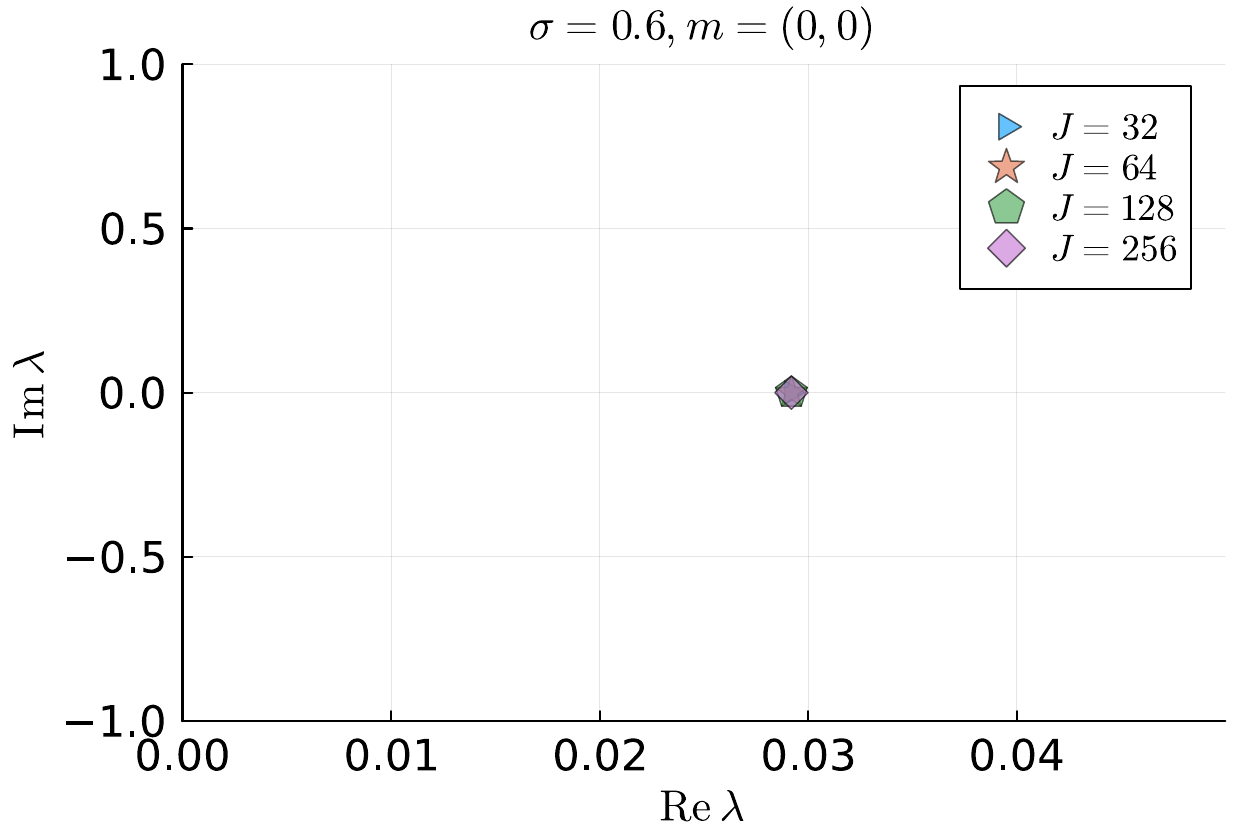}}
\subfigure[$\sigma=1.0$]{\includegraphics[width=4cm]{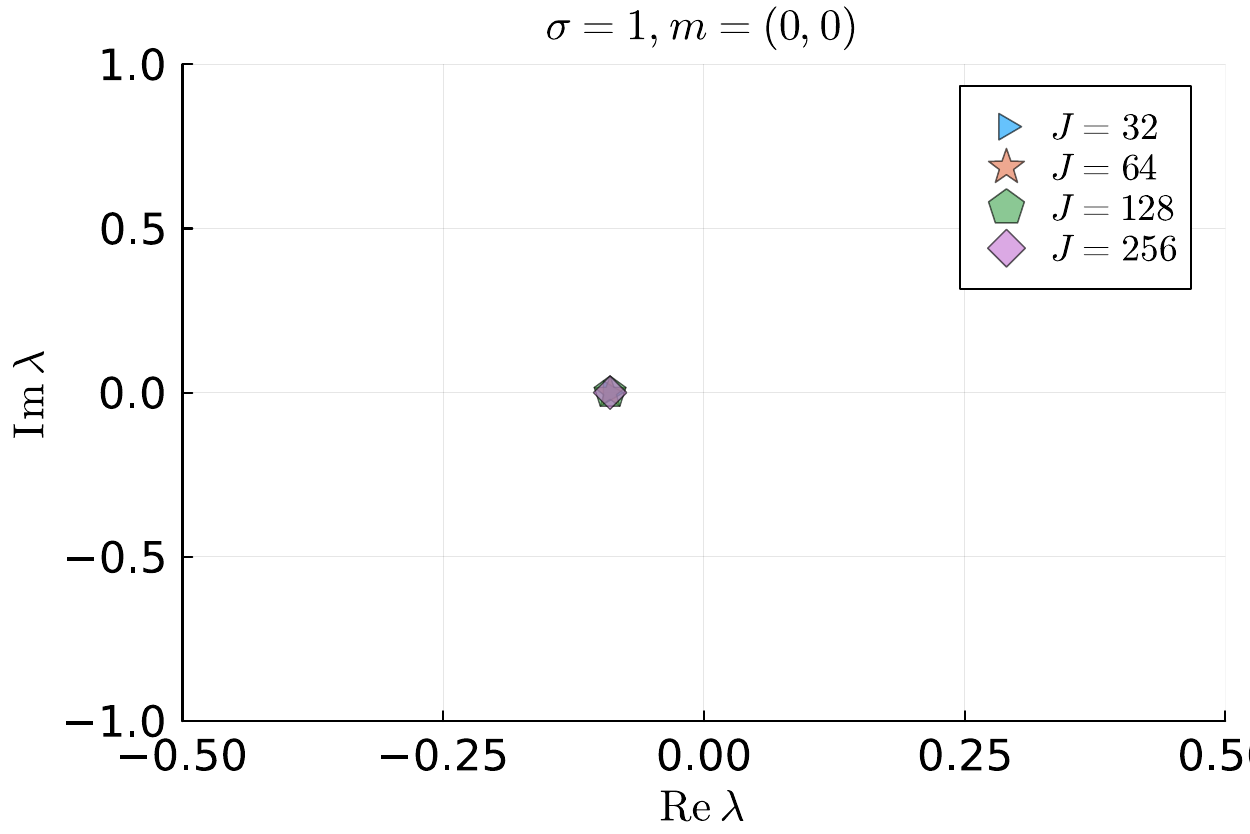}}
\caption{Eigenvalue of largest real part for the ${(m_1,m_2)}=(0,0)$ solution at different $\sigma$.  In the subcritical regime, this solution is linearly unstable, while in the supercritical regime it appears stable.}
\label{f:spec1}
\end{figure}

As we see in Figure \ref{f:spec1} - {whose numerical simulations have been carried out for different values of $\sigma$ with the choice $(m_1,m_2)=(0,0)$ i.e. $\rho_{\infty}(x)=\frac{1}{Z_{\sigma}}e^{-\frac{1}{\sigma}\cos\,2x}$} - there is an eigenvalue of positive real part for $\sigma$ in the subcritical regime.  This corresponds to a linear instability about this solution, and explains why we do not observe it in our time dependent simulations.  In contrast, for $\sigma$ in the supercritical regime, the eigenvalue nearest the imaginary axis has negative real part, and it is linearly stable; we do see this solution in the time dependent simulations.  Though the eigenvalues appear stable to numerical refinement, it remains to establish these properties for \eqref{linearised}.  We note that we filtered out zero eigenvalues of $\mathbf{L}$ as these corresponded to eigenstates that violated the mean zero condition of \eqref{linearised}.

\begin{figure}
\subfigure[$\sigma=0.2$]{\includegraphics[width=6cm]{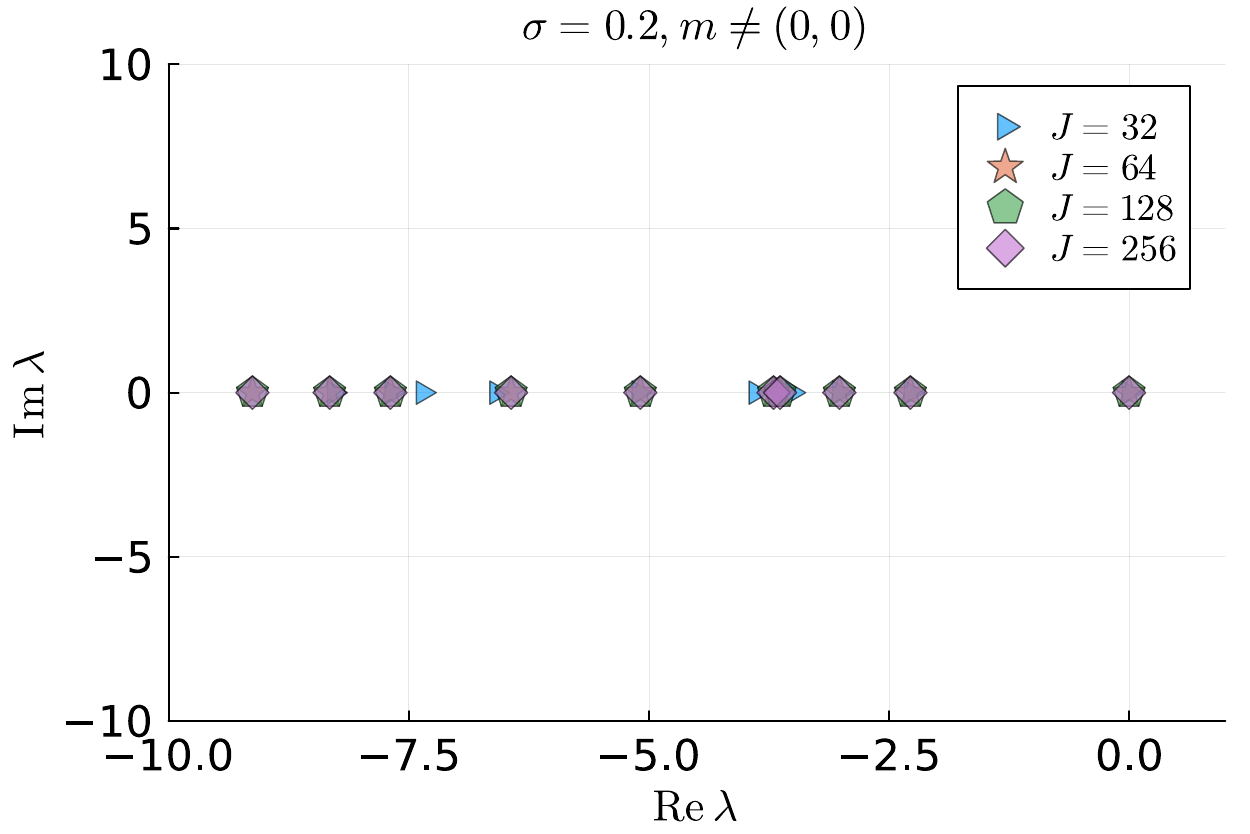}}
\subfigure[$\sigma=0.6$]{\includegraphics[width=6cm]{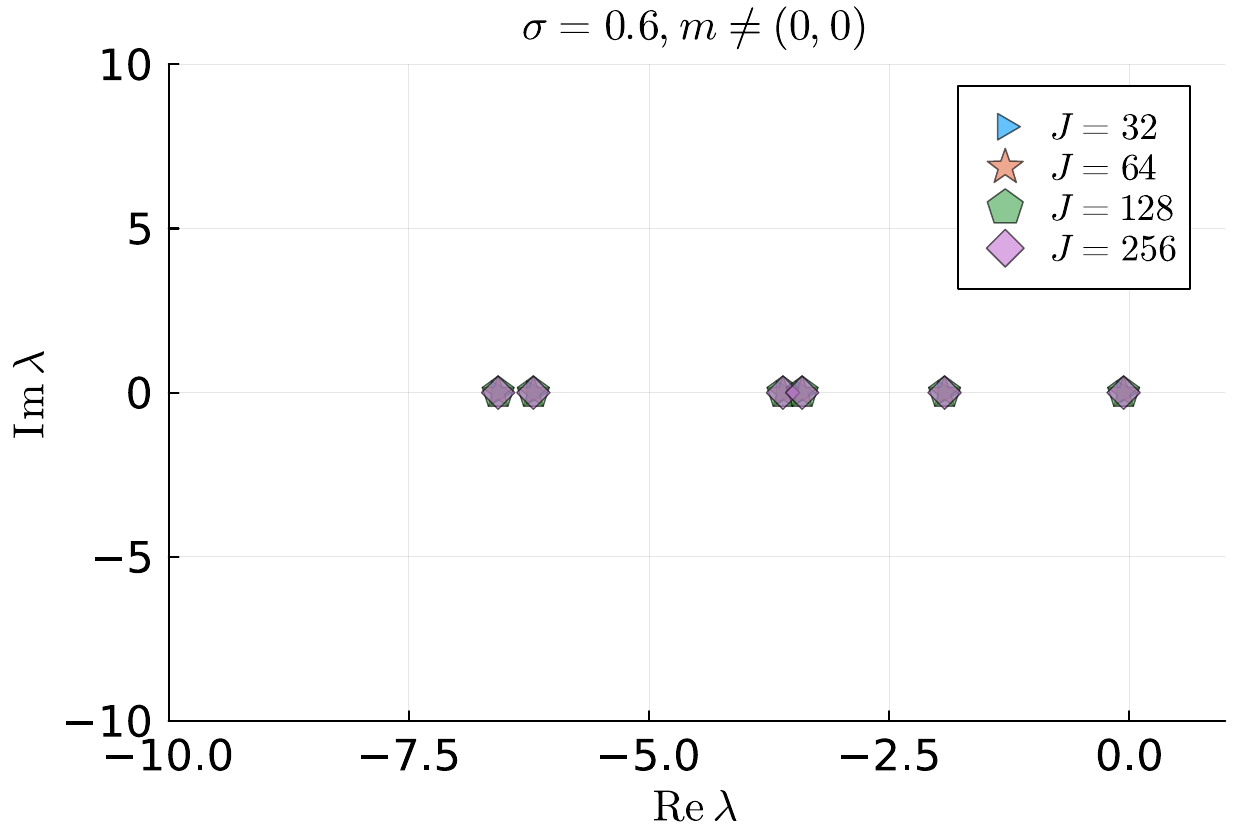}}

\subfigure[Zoom of $\sigma=0.2$]{\includegraphics[width=6cm]{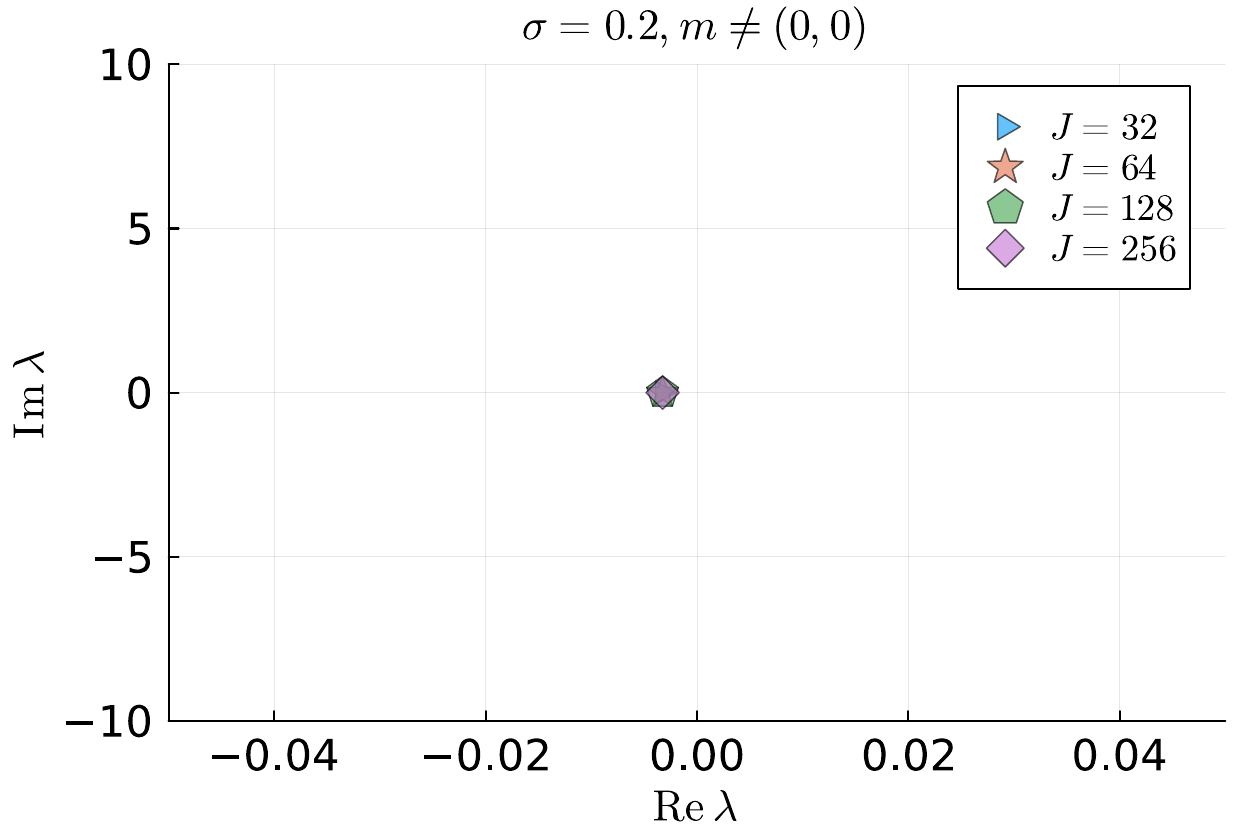}}
\subfigure[Zoom of $\sigma=0.6$]{\includegraphics[width=6cm]{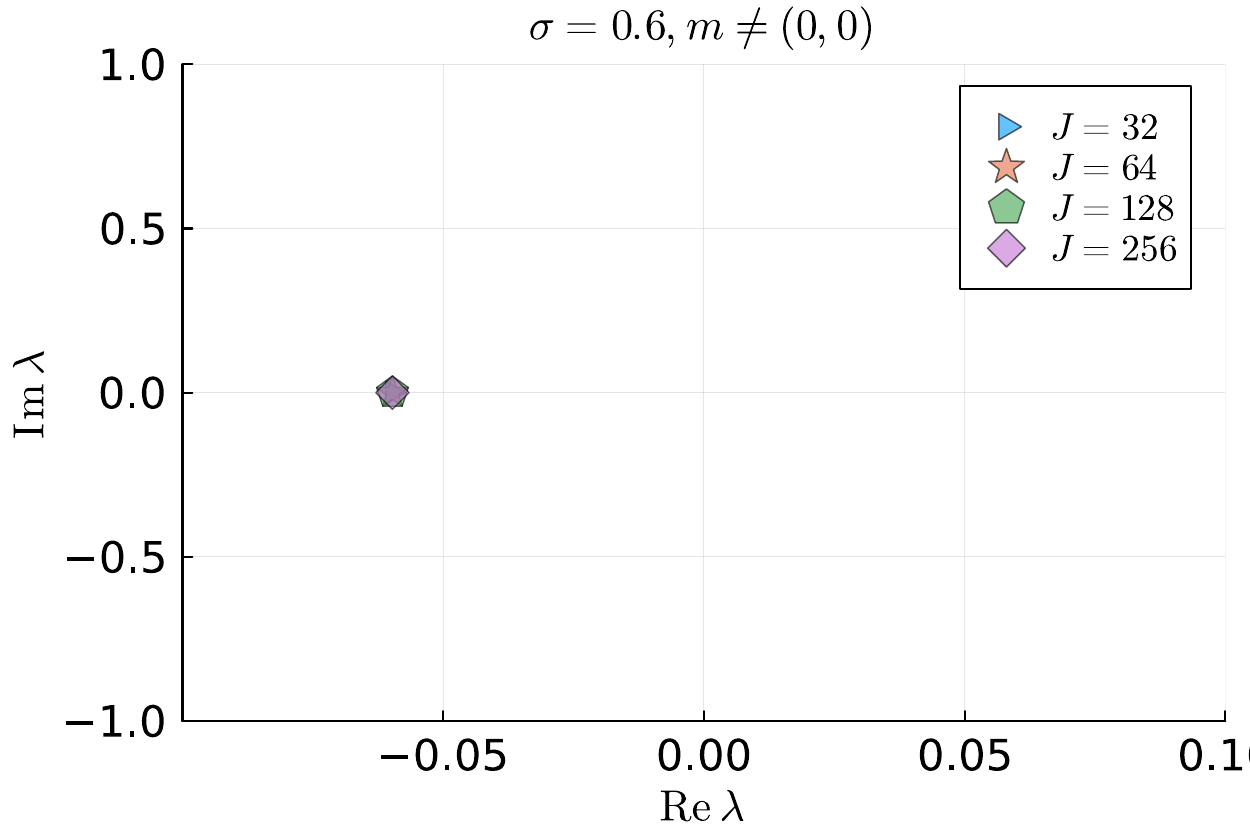}}

    \caption{Spectrum of the nontrivial stationary state, $m \neq (0,0)$, in the subcritical regime.  The eigenvalues nearest the real axis are negative, suggesting linear stability.}
    \label{fig:spec2}
\end{figure}

For the non-trivial solution, (i.e., ${(m_1,m_2)} \neq (0,0)$), in the subcritical regime, Figure \ref{fig:spec2} implies these are stable solutions.  As the images show, the eigenvalue nearest the imaginary axis is negative.  Thus, these will become metastable when stochastic forcing is added.

\section*{Acknowledgments.}

M.K. acknowledges support from the EPSRC grant EP/V520044/1 and support from the Italian Ministry of University and Research (MUR) via PRIN 2022– Project Title ConStRAINeD – CUP-2022XRWY7W. G.S. acknowledges support from the United States National Science Foundation,  grant number DMS-2111278. Work reported here was run on hardware supported by Drexel’s University Research Computing Facility. M.O. gratefully acknowledges support from the EPSRC grant EP/W034220/1 and from the Royal Society of Edinburgh Personal Research Fellowship scheme. 

\bibliographystyle{abbrv}
\bibliography{bib}
\end{document}